\documentclass[11pt]{article}
\usepackage{graphicx}
\usepackage{color}
\usepackage{amsmath}
\usepackage{amssymb}
\usepackage{amscd}
\usepackage{bbm}
\newcommand{\R}{\mathbb{R}}
\newcommand{\inr}[1]{\bigl< #1 \bigr>}

\newcommand{\E}{\mathbb{E}}

\newcommand{\eps}{\varepsilon}

\newtheorem{Theorem}{Theorem}[section]
\newtheorem{Lemma}[Theorem]{Lemma}

\newtheorem{Definition}[Theorem]{Definition}

\newtheorem{Corollary}[Theorem]{Corollary}

\numberwithin{equation}{section}

\def \proof {\noindent {\bf Proof.}\ \ }

\def \endproof
{{\mbox{}\nolinebreak\hfill\rule{2mm}{2mm}\par\medbreak}}

\begin{document}
\title{Dvoretzky type theorems for subgaussian coordinate projections}
\author{
Shahar Mendelson\footnote{Department of
Mathematics, Technion, I.I.T, Haifa 32000, Israel
\newline
{\sf email: shahar@tx.technion.ac.il}
\newline
Supported in part by the Mathematical Sciences Institute,
The Australian National University, Canberra, ACT 2601,
Australia. Additional support was given by the Israel Science Foundation grant 900/10.
 } }

\medskip
\maketitle
\begin{abstract}
Given a class of functions $F$ on a probability space $(\Omega,\mu)$, we study the structure of a typical coordinate
projection of the class, defined by $\{(f(X_i))_{i=1}^N : f \in
F\}$, where $X_1,...,X_N$ are independent, selected according to
$\mu$. This notion of projection generalizes the standard linear
random projection used in Asymptotic Geometric Analysis.

We show that when $F$ is a subgaussian class of functions, a typical coordinate projection satisfies a
Dvoretzky type theorem.
\end{abstract}

\section{Introduction}
Random projections appear naturally in various areas of mathematics, most notably in Asymptotic Geometric Analysis. What is arguably the most important result in classical Asymptotic Geometric Analysis, Milman's version of Dvoretzky's Theorem \cite{Mil71,MilSch,Pis}, deals with random sections/projections of a convex, centrally symmetric set in $\R^n$ with a nonempty interior (a convex body). The question was to identify the dimension $k$ for which an orthogonal projection of a convex body $T$ onto a typical element in the Granssmann manifold $G_{k,n}$, relative to the Haar measure, is almost Euclidean.

Milman showed that $k$ is governed by two parameters: the mean-width of $T$, defined by
$$
\int_{S^{n-1}} \sup_{t \in T} \inr{t,x}dx,
$$
where the integration is with respect to the Haar measure on the sphere; and the Euclidean radius of $T$, $\sup_{t \in T} \|t\|_{\ell_2^n}$, denoted by $d_T$.

An accurate formulation of Milman's Theorem for a gaussian projection (see, e.g.,
\cite{Pis}) is the following. Let $G=(g_i)_{i=1}^n$ be the standard gaussian vector on $\R^n$, whose coordinates are independent, standard gaussian random variables. Set
$$
\ell_*(T)=\E \sup_{t \in T} \sum_{i=1}^n g_i t_i,
$$
the gaussian mean-width of $T$. The critical dimension of $T$ is defined to be
$$
k_T^*=\left(\frac{\ell_*(T)}{d_T}\right)^2.
$$
Let $G_1,...,G_N$ be
independent copies of $G$ and put $\Gamma = \sum_{i=1}^N \inr{G_i,\cdot}e_i$. For every $0<\eps<1/2$, let $k_{T,\eps}^*=\eta_\eps k_T^*$, where
$\eta_\eps=\eps^2/\log(1/\eps)$. Finally, set $\ell_2^N$ to be $\R^N$ endowed with the Euclidean norm, and put $B_2^N$ to be the unit ball in $\ell_2^N$.  

\begin{Theorem} \label{thm:gaussian-Dvo-proj}
There exist absolute constants $c_1$ and $c_2$ for which the
following holds. If $0<\eps<1/2$ and $N=c_1k_{T,\eps}^*$, then
with probability at least $1-2\exp(-c_2k_{T,\eps}^*)$,
$$
(1-\eps)\ell_*(T) B_2^N \subset \Gamma T \subset (1+\eps)\ell_*(T)B_2^N.
$$
\end{Theorem}
Therefore, with high probability, a gaussian projection of $T$ of dimension proportional to $k_{T,\eps}^*$ is `almost' the Euclidean ball. In fact, the dependence on $\eps$ was improved further by Gordon \cite{Gor1,Gor2} to $c\eps^2k_T^*$, but since the focus of this note is on isomorphic results rather than on almost isometric ones, we will not go into more details on that dependence.

\vskip0.5cm

It is interesting to note that the proof of Theorem \ref{thm:gaussian-Dvo-proj} is indirect, and follows by dualizing the corresponding result for sections. Moreover, the proof is based rotation invariance and on a concentration theorem for Lipschitz functions on $\R^n$, relative to the standard gaussian measure.  Since the two are rather special properties that are rarely satisfied by more general matrix ensembles, it is not clear whether a Dvoretzky type theorem is true for linear projections selected according to other distributions.

\vskip0.5cm
Random projections appear naturally in other types of problems, and
the ones that motivated this work originated in Probability/Statitics.

Consider a class of functions $F$ defined a probability
space $(\Omega,\mu)$. If $X$ is distributed according to $\mu$, $X_1,...,X_N$ are independent copies of $X$ and
$\sigma=(X_i)_{i=1}^N$, the corresponding coordinate projection of $F$ is
\begin{equation} \label{eq:P-sigma-F}
P_\sigma F =\{ (f(X_i))_{i=1}^N : f \in F\} \subset \R^N.
\end{equation}
The name `random coordinate projection' may appear a little misleading when one is used to the linear setup. However, this notion seems to be the right generalization of a random linear projection. Indeed, let $\mu$ be a measure on $\R^n$ and consider the random matrix ensemble consisting of the matrices $\Gamma =\sum_{i=1}^N \inr{X_i,\cdot}e_i$, for $\sigma=(X_i)_{i=1}^N$ that is selected according to $\mu^N$. For a set $T \subset \R^n$, let $F_T=\{\inr{t,\cdot} : t \in T\}$ be the class of linear functionals on $\R^n$ associated with $T$. Clearly,
$$
P_\sigma F_T = \Gamma T
$$
that is, the linear projection of $T$, $\Gamma T$, is the corresponding coordinate projection of $F_T$ generated by the sample $\sigma$.

\vskip0.5cm

The key fact behind the results presented below, and which leads to a Dvoretzky type theorem for various coordinate projections, is due to Rudelson and Vershynin \cite{RudVer}. They proved a Dvoretzky type theorem for sections of a convex body (see its formulation for coordinate projections in Theorem \ref{thm:RudVer}), by showing that a body that contains the Euclidean unit ball $B_2^n$, has a coordinate section of the `right dimension' $m$ which is contained in an appropriate multiple of $B_1^m$, the unit ball in $\ell_1^m$. This immediately leads to an isomorphic Dvoretzky type theorem because a typical proportional section of $B_1^m$ is actually Euclidean (see, e.g. \cite{MilSch,Pis}).

Here, we will use the dual formulation of Rudelson and Vershynin's result to show that if $F$ is an $L$-subgaussian class of functions (defined below) and satisfies an additional regularity assumption, one may obtain a Dvoretzky type result for a typical coordinate projection $P_\sigma F$.

\begin{Definition} \label{def:L-subgaussian}
Given a function $f$ on the probability space $(\Omega,\mu)$, let
$$
\|f\|_{\psi_2} = \inf \left\{c>0: \E \exp(|f/c|^2) \leq 2 \right\}.
$$
A class of functions $F$ on $(\Omega,\mu)$ is
$L$-subgaussian if for every $f,h \in F \cup \{0\}$, $\|f-h\|_{\psi_2}
\leq L \|f-h\|_{L_2}$, where both norms are with respect to the underlying measure $\mu$.
\end{Definition}

Let $F \subset L_2$ be an $L$-subgaussian class that is convex and centrally symmetric and let $\{G_f : f \in F\}$ be the canonical gaussian process indexed by $F$, i.e., its covariance structure coincides with $L_2(\mu)$. In such a case, the natural analogs for $\ell_*$ and $d$ are $\E \sup_{f \in F} G_f \equiv \E \|G\|_F$ and $d_F = \sup_{f \in F} \|f\|_{L_2}$ respectively. Thus,
$$
k_F^* = \left(\frac{\E\|G\|_F}{d_F}\right)^2
$$
is the function-class analog of the critical dimension.

To formulate the main result of this note, consider a class $F$, a sample $(X_i)_{i=1}^N$ and a subset $I \subset \{1,...,N\}$. Let $Q_I : \R^N \to \R^I$ be defined for every $x=\sum_{i=1}^N x_ie_i$ by $Q_Ix=\sum_{i \in I} x_i e_i$. Set $V=P_\sigma F$ and $Q_I V = \{(f(X_i))_{i \in I} : f \in F\} \subset \R^I$. Let  $B_\infty^I$ be the unit cube on the coordinates $I$ and put $B_2^I$ to be the unit Euclidean ball on those coordinates.

Finally, let $D$ be the unit ball of $L_2(\mu)$ and set
$\phi(r) =\E \sup_{f \in F \cap r D} G_f$, the oscillation function of the gaussian process indexed by $F$.

\vskip0.5cm

\noindent{\bf Theorem A.}
For every $L \geq 1$ and $0<\alpha<1$ there exist constants $c_1,c_2,c_3,c_4$ and $c_5$ that depend only on $L$ and $\alpha$ for which the following holds. Let $F$ be an $L$-subgaussian class of functions that is convex and centrally symmetric. Assume further that $\phi(\alpha d_F) \leq \E\|G\|_F/4$. If $N \geq c_1 k^*_F$, then with $\mu^N$-probability at least $1-2\exp(-c_2k_F^*)$ there exists $I \subset \{1,...,N\}$, $|I| \geq c_3 k^*_F$ for which
\begin{equation*}
c_4\frac{\E\|G\|_F}{\sqrt{|I|}} B_\infty^I \subset Q_I V \subset c_5 \E\|G\|_F B_2^I.
\end{equation*}

\vskip0.5cm

The significant difference between a Dvoretzky type theorem and Theorem A is that the latter ensures the existence of an extremal cube contained in $Q_I V \subset \E\|G\|_F B_2^I$, rather than in the ball $c\E\|G\|_F B_2^I$ itself (clearly, $B_\infty^I/\sqrt{|I|}$ is the largest possible cube that one may find in $B_2^I$).

Just as noted above regarding the result from \cite{RudVer}, Theorem A is not far from  a Dvoretzky type theorem. Indeed, it is well known that a linear random projection (e.g. relative to the Haar measure -- but also with respect to more general random ensembles, as will be shown later) of the cube $B_{\infty}^I/\sqrt{|I|}$ is equivalent to the Euclidean ball $B_2^I$. Hence, Theorem A implies that $Q_I V$ is a subset of $\R^I$ that is a proportional random projection away from an isomorphic equivalence with a Euclidean ball.

\vskip0.5cm

The existence of extremal structures (in this case, of an extremal
cube) in a convex body, usually occurs because the set is, on one
hand, well bounded, and on the other, of extremal complexity. The
combination of the two properties forces some structure to appear.
Here, the rather weak assumption on the oscillation function $\phi(r)$ is used to ensure that for a typical $\sigma$,
\begin{equation} \label{eq:condition-general-D}
\ell_*(P_\sigma F) \geq c(L,\alpha) \sqrt{N} \E\|G\|_F
\end{equation}
and thus $P_\sigma F$ is a convex subset of $c\E\|G\|_F B_2^N$ of extremal gaussian width.

\vskip0.5cm

In addition to Theorem A, we will present two other applications when the class of functions is $F_T=\{\inr{t,\cdot} : t \in T\}$ for a convex body $T \subset \R^n$, and $\mu$ is an isotropic measure on $\R^n$ (recall that a probability measure $\mu$ on $\R^n$ is isotropic if it is symmetric and for every $t \in \R^n$, $\int_{\R^n} \inr{x,t}^2d\mu(x) = \|x\|_{\ell_2^n}^2$). The first application leads to a subgaussian Dvoretzky type theorem for spaces with a nontrivial cotype 2 constant; the second studies linear subgaussian images of the intersection body of the unit ball of $\ell_1^n$ with a Euclidean ball, and in particular, provides some information on the structure of certain random polytopes. Both results follow from appropriate versions of Theorem A, though not directly from Theorem A itself.

\vskip0.5cm

We end the introduction with a few basic definitions, some notation and facts that will be used throughout the note.
Absolute constants are denoted by $c_1,c_2,...$; their value may change from line to line. We write $A \lesssim B$ if there is an absolute constant $c_1$ for which $A \leq c_1 B$, and $A \sim B$ if $c_1 A \leq B \leq c_2 A$ for absolute constants $c_1$ and $c_2$.  $A \lesssim_r B$ or $A \sim_r B$ means that the constants depend on some parameter $r$.

Given a probability measure $\mu$ and $\alpha \geq 1$, $L_{\psi_\alpha}$ is the Orlicz space of all measurable functions, for which the $\psi_\alpha$ norm, defined by
$$
\|f\|_{\psi_\alpha}=\inf\left\{ c>0 : \E_\mu \exp(|f/c|^\alpha) \leq 2 \right\},
$$
is finite. Basic facts
on Orlicz spaces may be found in \cite{VW}.

One feature of a $\psi_\alpha$ random variable is that an average of its independent copies concentrates around its mean.

\begin{Theorem} \label{thm:basic-probab-inequality}
There exists an absolute constant $c_1$ for which the following holds. If $f \in L_{\psi_1}$ and $X_1,...,X_N$ are independent random variables distributed according to $\mu$, then for every $u>0$,
$$
Pr\left(\left|\frac{1}{N} \sum_{i=1}^N f(X_i) -\E f\right| \geq u \|f\|_{\psi_1} \right) \leq 2 \exp(-c_1N \min\{u^2,u\}).
$$
In particular, if $f \in L_{\psi_2}$ then
$$
Pr\left(\left|\frac{1}{N} \sum_{i=1}^N f^2(X_i) -\E f^2\right| \geq u \|f\|^2_{\psi_2} \right) \leq 2 \exp(-c_1N \min\{u^2,u\}).
$$
\end{Theorem}
The first part of Theorem \ref{thm:basic-probab-inequality} is a $\psi_1$ version of Bernstein's inequality (see, for example, \cite{VW}); the second one is an immediate outcome of the first, because $\|f^2\|_{\psi_1}=\|f\|_{\psi_2}^2$.

Finally, if $(a_i)_{i=1}^N \in \R^N$, denote by $(a_i^*)_{i=1}^N$ a monotone non-increasing rearrangement of $(|a_i|)_{i=1}^N$.

\section{Remarks on a Dvoretzky type theorem for gaussian projections} \label{sec:gaussian}
In this section we will sketch the argument behind Milman's version of Dvoretzky's Theorem for gaussian projections. All the facts presented here are known, and will only serve as an indication of how a Dvoretzky type theorem may be extended to the case we are interested in: typical coordinate projections of a function class.

For reasons that will become clear later, the argument will be split into
two parts. The first is the upper estimate that follows from information on the monotone rearrangement of the
random vectors $(\inr{G_i,t})_{i=1}^N$:

\begin{Theorem} \label{lemma:monotone-gaussian}
There exists absolute constants $c_1,c_2$ and $c_3$ for which the
following holds. For every $\eps>0$, $1 \leq k \leq N$ and $u \geq
c_1$, with probability at least $1-2\exp(-c_2u^2 k\log (eN/k
\eps))$,
$$
\sup_{t \in T} \left(\sum_{i=1}^k (\inr{G_i,t}^*)^2  \right)^{1/2}
\leq (1+\eps) \left(\ell_*(T) + c_3ud_T \sqrt{k
\log(eN/k\eps)}\right).
$$
In particular, with probability at least $1-2\exp(-c_2u^2N)$,
\begin{equation} \label{eq:gaussian-diameter}
\sup_{t \in T} \left(\sum_{i=1}^N \inr{G_i,t}^2 \right)^{1/2} \leq
(1+\eps)\left(\ell_*(T) + c_3ud_T \sqrt{N}\right).
\end{equation}
\end{Theorem}

The upper estimate in Dvoretzky's Theorem follows from Theorem \ref{lemma:monotone-gaussian}. Indeed, if $\Gamma=\sum_{i=1}^N \inr{G_i,\cdot}e_i$ then
\begin{align*}
\sup_{t \in T} \|\Gamma t\|_{\ell_2^N} \leq & \left(1+\eps\right) \cdot \left(\ell_*(T)+c_1ud_T \sqrt{N\log(e/\eps)} \right)
\\
\leq &(1+\eps)\ell_*(T) \left(1+c_1v \frac{d_T}{\ell_*(T)}\sqrt{N\log(e/\eps)}\right)
\\
\leq & (1+\eps)\ell_*(T) \cdot(1+c_1v\eps),
\end{align*}
provided that $N \lesssim k_{T,\eps}^*$.

Note that a proof of an isomorphic upper estimate is an immediate outcome of   \eqref{eq:gaussian-diameter}. Hence, a high probability estimate of the form
$$
\sup_{f \in F} \left(\sum_{i=1}^N f^2(X_i) \right)^{1/2} \leq c\left(\E \|G\|_F + d_F\sqrt{N}\right)
$$
implies that if $|\sigma| \leq c_1k_F^*$, then $P_\sigma F \subset c\E\|G\|_F B_2^{|\sigma|}$ for a typical projection.

\vskip0.5cm

The other half of Theorem \ref{thm:gaussian-Dvo-proj} turns out to be more restrictive.
\begin{Theorem} \label{thm:lower-D-gaussian}
There exist absolute constants $c_1$ and $c_2$ for which the following holds. If $0<\eps<1/2$, $N \leq c_1 k_{T,\eps}^*$ and $\Gamma=\sum_{i=1}^N \inr{G_i,\cdot}e_i$, then with probability at least
$1-2\exp(-c_2\eps^2k_{T,\eps}^*)$,
$$
(1-\eps) \ell_*(T) B_2^N \subset \Gamma T.
$$
\end{Theorem}

The proof is based on a separation argument:
Fix $\rho>0$ and an integer $N$. If $\rho B_2^N \not \subset \Gamma T$, there is $w \in \rho B_2^N \backslash \Gamma T$, and since $\Gamma T$ is a convex body, there is a functional $z \in S^{N-1}$ for which $\sup_{t \in T} \inr{\Gamma t,z} < \inr{w,z}$. Clearly, $\inr{w,z} \leq \rho$, and thus it suffices to show that for the right choice of $N$, with high probability,
\begin{equation} \label{eq:Dvo-lowe-in-proof}
\inf_{z \in S^{N-1}} \sup_{t \in T} \inr{\Gamma t,z} = \inf_{z \in S^{N-1}} \|\sum_{i=1}^N z_i G_i \|_{T^\circ} > \rho,
\end{equation}
where $\| \ \|_{T^\circ}$ is the norm on $\R^n$ whose unit ball is $T^\circ$, the polar body of $T$.

\vskip0.5cm

Observe that \eqref{eq:Dvo-lowe-in-proof} actually follows from a small-ball estimate rather than a concentration based one. It holds for $\rho \sim \ell_*(T)$ and the right choice of $N$, if for every $z \in
S^{N-1}$ and $u<1/2$,
$$
Pr\left(\sup_{t \in T} \left|\sum_{i=1}^N z_i \inr{G_i,t}\right| \leq cu\ell_*(T)
\right) \leq u^{k_T^*}.
$$
Although a small-ball estimate of this type is not unique to the gaussian
ensemble, it is still rather restrictive, certainly in the context of coordinate projections of function classes. Building on the result from \cite{RudVer}, we will explain why, instead of a small-ball condition, and once the `upper estimate' in Theorem A is satisfied, the `lower estimate' holds when
$$
\int_{S^{N-1}} \sup_{f \in F} \left|\sum_{i=1}^N z_i f(X_i)\right| dz \geq c\E\|G\|_F.
$$

\section{A few facts on chaining} \label{sec:pre}
Our results are based on chaining methods and we refer the reader to \cite{Tal:book} for an extensive survey on this topic.

\begin{Definition} \label{def:gamma-2} \cite{Tal:book}
For a metric space $(F,d)$, an {\it admissible sequence} of $F$ is a
collection of subsets of $F$, $\{F_s : s \geq 0\}$, satisfying that for
every $s \geq 1$, $|F_s| \leq 2^{2^s}$ and $|F_0|=1$. For $s_0 \geq 0$, let
$$
\gamma_{2,s_0}(F,d) =\inf \sup_{f \in F} \sum_{s=s_0}^\infty
2^{s/2}d(f,F_s),
$$
where the infimum is taken with respect to all admissible sequences
of $F$.

If $s_0=0$ we will write $\gamma_2(F,d)$ instead of $\gamma_{2,s_0}(F,d)$.
\end{Definition}
If $F$ is a class of functions and $(F_s)_{s \geq 0}$ is an admissible sequence, let $\pi_s f$ be a nearest point to $f$ in $F_s$ relative to the metric $d$, and for $s>0$, let $\Delta_s f = \pi_s f - \pi_{s-1}f$.

\vskip0.5cm

When $F \subset L_2(\mu)$ (or $\ell_2^N$), $\gamma_2(F,L_2)$ is determined by properties of the canonical gaussian process indexed by the class (see \cite{Dud:book,Tal:book} for detailed
expositions on these connections). Indeed, under certain mild measurability
assumptions, if $\{G_f: f \in F\}$ is a centered gaussian process indexed by $F$, then setting $\E\|G\|_F \equiv \E \sup_{f \in F} G_f$ one has
$$
c_1 \gamma_2(F,d) \leq \E\|G\|_F \leq c_2 \gamma_2(F,d),
$$
where $c_1$ and  $c_2$ are  absolute constants, and for every $f,h
\in F$, $d^2 (f,h) = \E|G_f-G_h|^2$.  The upper bound is due to
Fernique \cite{F} and the lower bound is Talagrand's Majorizing
Measures Theorem \cite{Tal87,Tal:book}. Hence, if $\{G_f : f \in F\}$ is the canonical gaussian process indexed by $F \subset L_2(\mu)$ then $\E \|G\|_F \sim \gamma_2(F,L_2)$.

Note that if $T \subset \ell_2^N$,
$(g_i)_{i=1}^N$ are independent, standard gaussian random variables and $G_t =
\sum_{i=1}^n g_i t_i$, then $d(s, t) = \|s-t\|_{\ell_2^n}$; therefore
\begin{equation}
  \label{maj_meas}
  c_1 \gamma_2(T,\ell_2^n) \leq  \E \sup_{t \in T} \sum_{i=1}^n
  g_it_i  \leq c_2 \gamma_2(T,\ell_2^n).
\end{equation}
Also, if $\mu$ is an isotropic measure on $\R^n$, $T \subset \R^n$ and $F_T=\{\inr{t,\cdot} : t \in T\} \subset L_2(\mu)$, the canonical gaussian process indexed by $F_T$ satisfies $\ell_*(T)=\E\|G\|_{F_T} \sim \gamma_2(T,\ell_2^n)$.

\section{A subgaussian Dvoretzky type theorem} \label{sec:beyond-gaussian-1}
As mentioned earlier, our results are based on \cite{RudVer}. Although not stated in exactly this way in \cite{RudVer}, Theorem \ref{thm:RudVer} follows immediately from Theorem 7.4 and Corollary 7.9 there:

\begin{Theorem} \label{thm:RudVer}
There exist absolute constants $c_1$ and $c_2$ for which the following holds. Let $V \subset \R^N$ and assume that
\begin{description}
\item{1.} $d_V \leq \alpha$, and
\item{2.} $\ell_*(V) \geq \delta \sqrt{N} \alpha$.
\end{description}
Then, there exists $I \subset \{1,...,N\}$, for which $|I| \geq c_1\frac{\delta^2N}{\log^3 (2/\delta)}$ and
$$
c_2 \frac{\alpha}{\sqrt{N}} B_\infty^I \subset Q_I V.
$$
\end{Theorem}

The sets $V$ we will be interested in are the random coordinate projections $P_\sigma F$, which leads to the following definition:
\begin{Definition} \label{def:class-A}
For every $u>0$, $0<\delta<1$ and a fixed integer $N$, let
${\cal A}_{u,\delta,N} \subset \Omega^N$ be the event on which
\begin{description}
\item{1.} For every $f \in F$,
$\left(\sum_{i=1}^N f^2(X_i)\right)^{1/2} \leq u\left(\E \|G\|_F + d_F \sqrt{N}\right)$,
and
\item{2.} $\ell_*(P_\sigma F) \geq \delta \sqrt{N}\E \|G\|_F$.
\end{description}
\end{Definition}
To simplify notation we will sometimes omit the subscripts $u,\delta$ and $N$.

\vskip0.5cm

The following is a direct outcome of Theorem \ref{thm:RudVer}.
\begin{Theorem} \label{thm:beyond-gaussian-main}
For every $0<\delta<1$ and $u>0$, there exist constants $c_1$, $c_2,c_3$ and $c_4$ that depend only on $\delta$ and $u$ for which the following holds. If $N \geq c_1 k^*_F$, $\sigma \in {\cal A}_{u,\delta,N}$ and $V = P_\sigma F$, then there exists $I \subset \{1,...,N\}$, $|I| \geq c_2 k^*_F$ for which
\begin{equation} \label{eq:beyond-gaussian-full}
c_3\frac{\E\|G\|_F}{\sqrt{|I|}} B_\infty^I \subset Q_I V \subset c_4 \E\|G\|_F B_2^I.
\end{equation}
In particular, if $\mu$ is isotropic, $T \subset \R^n$ and $\Gamma_I=\sum_{i \in I} \inr{X_i,\cdot}e_i$, then
\begin{equation} \label{eq:beyond-gaussian-R-n}
c_3\frac{\ell_*(T)}{\sqrt{|I|}} B_\infty^I \subset \Gamma_I T \subset c_4 \ell_*(T) B_2^I.
\end{equation}
\end{Theorem}

Thus, one has to identify conditions in which the event ${\cal A}$ has sufficiently high probability.

\vskip0.5cm
\subsection{The upper estimate}

\begin{Lemma} \label{lemma:subgaussian-diameter}
For every $L \geq 1$ there exist constants $c_1$ and $c_2$ that depend only on $L$ for which the following holds. If $F$ is an $L$-subgaussian class, then for every $u \geq 1$, with probability at least $1-2\exp(-c_1u^2 N)$,
$$
\sup_{f \in F} \left(\sum_{i=1}^N f^2(X_i) \right)^{1/2} \leq c_2u\left(\E\|G\|_F + d_F \sqrt{N}\right).
$$
\end{Lemma}
\proof Let $(F_s)_{s \geq 0}$ be an admissible sequence of $F$ and fix $s_0$ to be the first integer $s$ for which $2^s \geq N$. Given a sample $X_1,...,X_N$, set $\|f\|_{L_2^N} = (N^{-1}\sum_{i=1}^N f^2(X_i))^{1/2}$. Recall that $\Delta_s f = \pi_s f - \pi_{s-1} f$ and note that $f = \sum_{s>s_0} \Delta_s f + \pi_{s_0} f $. Therefore,
$$
\|f\|_{L_2^N} \leq \sum_{s>s_0} \|\Delta_s f \|_{L_2^N} + \|\pi_{s_0} f\|_{L_2^N}.
$$
Observe that if $h \in L_{\psi_2}$, then by Theorem \ref{thm:basic-probab-inequality}, with probability at least $1-2\exp(-cN\min\{v,v^2\})$,
$$
\frac{1}{N}\sum_{i=1}^N h^2(X_i) \leq \E h^2 + v \|h\|_{\psi_2}^2.
$$
Let $u \geq 1$, set $v=u 2^{s}/N \geq 1$ and put $h=\Delta_s f$ (resp. $h = \pi_{s_0} f$). Since $F$ is $L$-subgaussian, $\|\Delta_s f\|_{\psi_2} \leq L\|\Delta_s f\|_{L_2}$  and $\|\pi_{s_0} f\|_{\psi_2} \leq L\|\pi_{s_0} f\|_{L_2}$. Moreover, as $|F_{s-1}| \cdot |F_{s}| \leq 2^{2^{s+1}}$, if $u \geq c_0$ then with probability at least $1-2\exp(-c_1u 2^{s_0}) \geq 1-2\exp(-c_2uN)$, for every $f \in F$ and every $s>s_0$,
$$
\sum_{i=1}^N (\Delta_s f)^2(X_i) \leq (N+L^2u2^s)\|\Delta_s f\|_{L_2}^2,
$$
and
$$
\sum_{i=1}^N (\pi_{s_0} f)^2(X_i) \leq (N+L^2u2^{s_0})\|\pi_{s_0} f\|_{L_2}^2.
$$
On that event,
\begin{align*}
& \sum_{s > s_0} \|\Delta_s f\|_{L_2^N}+\|\pi_{s_0} f\|_{L_2^N} \lesssim L\sqrt{u}(\sum_{s>s_0} 2^{s/2}\|\Delta_s f\|_{L_2} + 2^{s_0/2}d_F)
\\
\lesssim & L\sqrt{u} \left(\gamma_{2,s_0}(F,L_2) + 2^{s_0/2}d_F\right) \lesssim L\sqrt{u} \left(\E\|G\|_F + \sqrt{N}d_F\right),
\end{align*}
provided that $(F_s)_{s \geq 0}$ is an almost optimal admissible sequence with respect to the $L_2$ norm.
\endproof
\vskip0.5cm

\subsection{The lower estimate}

Next, we turn to the second, more restrictive condition in the definition of ${\cal A}$.

\subsubsection*{Classes with a well-behaved gaussian oscillation}
Let $F \subset L_2(\mu)$ be a convex and centrally symmetric class. Recall that $\phi(r)=\E \sup_{f \in F \cap rD} G_f$ and assume that there is $0<\alpha<1$ for which
$$
\phi(\alpha d_F) \leq \frac{\E\|G\|_F}{4}.
$$
Clearly, such an $\alpha$ exists if $\{G_f : f \in F\}$ is a continuous gaussian process.

The first observation needed for the proof of Theorem A is a standard subgaussian version of the Johnson-Lindenstrauss Lemma.
\begin{Lemma} \label{Lemma:extended-JL}
For every $L > 1$ there exist constants $c_1, c_2, c_3$ and $c_4$
that depend only on $L$ and for which the
following holds. If $H$ is an $L$-subgaussian class of functions
with $|H| \leq \exp(k)$, then
for every $N \geq c_1 k$, with $\mu^N$-probability at least $1-2\exp(-c_2N)$, for every $h_1,h_2 \in H$,
$$
\frac{1}{2}\|h_1-h_2\|_{L_2}^2 \leq \frac{1}{N}\sum_{i=1}^N (h_1-h_2)^2(X_i) \leq \frac{3}{2}\|h_1-h_2\|_{L_2}^2.
$$
In particular, on the same event,
$$
c_3 \sqrt{N}\E\|G\|_H \leq \ell_*(P_\sigma H) \leq c_4 \sqrt{N}\E\|G\|_H.
$$
\end{Lemma}

The first part follows from Theorem \ref{thm:basic-probab-inequality}, while the second is a corollary of the first part and Slepian's Lemma (see, e.g. \cite{LT}).

\vskip0.5cm
\noindent{\bf Proof of Theorem A.}
Let $\Lambda$ be a maximal $\alpha d_F$-separated subset of $F$ with respect to the $L_2(\mu)$ norm. Since $F$ is convex and centrally symmetric, for every $f \in F$ one has that  $f=\pi f + (f -\pi f)$ where $\pi f \in \Lambda$ and $f-\pi f \in 2F \cap \alpha d_F D \subset 2(F \cap \alpha d_f D)$. Therefore,
\begin{align*}
\E \|G\|_F \leq & \E \sup_{f \in
\Lambda} G_f + 2 \E \sup_{f \in F \cap \alpha d_F D} G_f = \E \sup_{f \in
\Lambda} G_f +2 \phi(\alpha d_F)
\\
\leq & \E \sup_{f \in \Lambda} G_f + \frac{1}{2} \E \|G\|_F,
\end{align*}
and thus $\E \sup_{f \in \Lambda} G_f \geq (1/2)\E \|G\|_F$. On
the other hand, by Sudakov's minoration (see, e.g., \cite{LT}),
$$
\log |\Lambda|
\leq c_0\left(\frac{\E \|G\|_F}{\alpha d_F}\right)^2 = (c_0/\alpha^2)k_F^*.
$$

Let $c_1,...,c_4$ as in Lemma \ref{Lemma:extended-JL} and note that by that lemma, applied to the set $\Lambda$ for $k=(c_0/\alpha^2) k_F^*$ and $N=c_1 k$, one has that with probability at least $1-2\exp(-c_2 N)$,
$$
c_3 \sqrt{N}\E \sup_{f \in
\Lambda} G_f \leq \ell_*(P_\sigma \Lambda) \leq c_4 \sqrt{N}\E \sup_{f \in
\Lambda} G_f.
$$
Hence, with probability at least $1-2\exp(-c_5(L,\alpha)k_F^*)$,
$$
\ell_*(P_\sigma F) \geq c_6(L,\alpha) \sqrt{N} \E\|G\|_F,
$$
implying that for $\delta \sim_{L,\alpha} 1$, $N \sim_{L,\alpha} k_F^*$ and $u=1$,
$$
Pr({\cal A}_{u,\delta,N}) \geq 1-2\exp(-c_7(L,\alpha) k_F^*).
$$
Therefore, by Theorem \ref{thm:beyond-gaussian-main}, if $\sigma \in {\cal A}$ and $V=P_\sigma F$, there is a subset $I \subset \{1,...,N\}$, $|I| \sim_L k_F^*$ and
$$
c_{8} B_\infty^I \frac{\E\|G\|_{F}}{\sqrt{k_F^*}} \subset Q_I V \subset c_{9} \E\|G\|_{F} B_2^I
$$
for constants $c_8$ and $c_9$ that depend only on $L$ and $\alpha$, as claimed.
\endproof

\subsubsection*{Spaces with cotype}
Let $T \subset \R^n$ be a convex body and assume that $\| \ \|_{T^\circ}$, the norm whose unit ball is the polar body $T^\circ$, has (gaussian) cotype 2 with a constant $C_2(T^\circ)$.

It is well known that if $X$ is an isotropic, $L$-subgaussian vector on $\R^n$ with iid coordinates, there exist constants $c_0=c_0(L)$ and $c_1$ that depend only on $C_2(T^\circ)$ for which
\begin{equation} \label{eq:cotype}
c_0 \E\|X\|_{T^\circ} \leq \ell_*(T) \leq c_1 \E \|X\|_{T^\circ}.
\end{equation}

The left-hand side of \eqref{eq:cotype} follows, for example, from a chaining argument and the Majorizing Measures Theorem, while the right-hand side may be found in \cite{LT} (see also \cite{Men-TJ}).

\begin{Lemma} \label{lemma:basic-estimate-cotype}
For every $L > 1$ and $\kappa>0$ there exists constants $c_1$ and $c_2$ that depend only on $L$ and $\kappa$ for which the following holds. If $T$ is a convex body for which $C_2(T^\circ) \leq \kappa$, then with $\mu^N$-probability at least $1-2\exp(-c_1N)$,
$$
\ell_*(P_\sigma T) \geq c_2\sqrt{N} \ell_*(T).
$$
\end{Lemma}

\proof For every $X_1,...,X_N$,
$$
\ell_*(P_\sigma T) = \E_g \sup_{t \in T} \sum_{i=1}^N g_i\inr{X_i,t}= \E_g\|\sum_{i=1}^N g_i X_i\|_{T^\circ}.
$$
By the Kahane-Khintchine inequality and since $T^\circ$ has cotype $2$, $$
\left(\E_g \|\sum_{i=1}^N g_i X_i\|_{T^\circ} \right)^2 \gtrsim  \E_g\|\sum_{i=1}^N g_i X_i\|^2_{T^\circ} \geq C_2^{-2}(T^\circ) \sum_{i=1}^N \|X_i\|_{T^\circ}^2.
$$
To conclude the proof one has to find a high probability lower bound on $\sum_{i=1}^N \|X_i\|_{T^\circ}^2$. To that end, observe that there are constants $c_1$ and $0<\eta<1$ that depend only on $L$ and $C_2(T^\circ)$, for which
\begin{equation} \label{eq:small-ball-cotype}
Pr\left(\|X\|_{T^\circ} \geq c_1\E \|X\|_{T^\circ} \right) \geq \eta.
\end{equation}
Indeed, let $(T_s)_{s \geq 0}$ be an almost optimal admissible sequence  of $T$ and set $s_0 \geq 0$. By a straightforward chaining argument one has that with probability at least $1-2\exp(-c_2u^22^{s_0})$, for every $t \in T$,
$$
|\inr{X,t}| \leq u \left(\sum_{s > s_0} 2^{s/2}\|\inr{\Delta_s t,X}\|_{\psi_2} + 2^{s_0/2}\|\inr{\pi_{s_0} t,X}\|_{\psi_2} \right).
$$
Since $X$ is $L$-subgaussian and by the Majorizing Measures Theorem, for every $u \geq 1$,
\begin{equation} \label{eq:tail-estimates-p}
Pr \left(\sup_{t \in T} |\inr{X,t}| \geq c_3Lu\left(\ell_*(T)+2^{s_0/2}d_T\right)\right) \leq 2\exp(-c_2u^2 2^{s_0}).
\end{equation}
If $p \geq 1$, set $2^{s_0} \sim \max\{p,k_T^*\}$. Integrating the tail estimate \eqref{eq:tail-estimates-p} implies that
$$
\left(\E\|X\|_{T^\circ}^p\right)^{1/p} \lesssim_L \ell_*(T) + \sqrt{p}d_T \leq c_4(\kappa,L) \left( \E \|X\|_{T^\circ} + \sqrt{p}d_T\right),
$$
and \eqref{eq:small-ball-cotype} follows from the Paley-Zygmund inequality.

Finally, if $(\eta_i)_{i=1}^N$ are selectors with mean $\eta$ (defined in \eqref{eq:small-ball-cotype}), then
\begin{align*}
& Pr \left(\sum_{i=1}^N \|X_i\|_{T^\circ}^2 \geq c_1^2(\E \|X\|_{T^\circ})^2N/100\right) \geq
Pr\left(\sum_{i=1}^N \eta_i \geq N\eta/100\right)
\\
\geq & 1-2\exp(-c_4N\eta^2).
\end{align*}
\endproof

Combining Lemma \ref{lemma:subgaussian-diameter} with Lemma \ref{lemma:basic-estimate-cotype} shows that for the
correct choice of $u$ and $\delta$, which depend only on $L$ and
on the cotype-2 constant of $T^\circ$, and for $N \sim_{L,C_2} k_T^*$,
$$
\mu^N({\cal A}_{u,\delta,N}) \geq 1-2\exp(-ck_T^*).
$$
And, if $\sigma \in {\cal A}$, $I$ as in Theorem \ref{thm:beyond-gaussian-main} and $\Gamma_I=\sum_{i \in I} \inr{X_i,\cdot}e_i$, one has
$$
c_1\frac{\ell_*(T)}{\sqrt{|I|}}B_\infty^I \subset \Gamma_I T \subset c_2\ell_*(T) B_2^I.
$$

\subsubsection*{$B_1^n$ and random polytopes}

Let $T=B_1^n$ be the unit ball in $\ell_1^n$ and recall the well known fact (see, for example, \cite{GLMP}) that for every $1/\sqrt{n} \leq \rho \leq
1$,
$$
\E \|G\|_{B_1^n \cap \rho B_2^n} =\E \sup_{t \in B_1^n \cap \rho B_2^n} \sum_{i=1}^n g_i t_i \sim
\sqrt{\log(en\rho^2)},
$$
while for
$\rho \lesssim 1/\sqrt{n}$, $\E \|G\|_{B_1^n \cap \rho B_2^n}
\sim \sqrt{n} \rho$. Therefore, if $\log n \lesssim k \leq n$, and $\rho_k
\sim \sqrt{\frac{\log(en/k)}{k}}$, the critical dimension of the intersection body $B_1^n \cap \rho_k B_2^n$ is
$$
\left(\frac{\ell_*(B_1^n \cap \rho_k B_2^n)}{\rho_k}\right)^2 = k.
$$

For every $I \subset \{1,...,n\}$, let $S^I=\{x \in S^{n-1} : {\rm supp}(x)=I\}$. Note that $\bigcup_{I} c_1 \rho_k
S^I \subset B_1^n \cap \rho_k B_2^n$, with the union taken over all
subsets of $\{1,...,n\}$ of cardinality $m \sim k/\log(en/k)$ and $c_1$ is an appropriate absolute constant. A standard argument (see, e.g., \cite{MPR}, Lemma 3.6) shows that there is a collection ${\cal B}$ of subsets of $\{1,...,n\}$ of cardinality $m$ that is $c_2m$-separated in the Hamming distance, and $\log |{\cal B}| \gtrsim k$. Let
$$
W_k=\left\{\frac{c_1\rho_k}{\sqrt{m}} \sum_{i \in I} e_i : I \in {\cal B}\right\} \subset \bigcup_I c_1 \rho_k S^I \subset B_1^n \cap \rho_k B_2^n,
$$
note that $|W_k| \leq \exp(c_3k)$ and that by Slepian's Lemma
$$
\E \sup_{w \in W_k} \sum_{i=1}^n g_i w_i \gtrsim \rho_k \sqrt{k} \gtrsim \ell_*(B_1^n \cap \rho_k B_2^n).
$$
Let $X$ be an isotropic, $L$-subgaussian vector on $\R^n$, distributed according to a measure $\mu$. Applying Lemma \ref{Lemma:extended-JL} to the set $W_k \subset B_1^n \cap \rho_k B_2^n$ for $N \sim k$ it follows that with $\mu^N$-probability at least $1-2\exp(-c_3(L)k)$,
\begin{equation} \label{eq:l-star-local-b-1-n}
\ell_*(P_\sigma(B_1^n \cap \rho_k B_2^n)) \gtrsim \sqrt{k} \ell_*(B_1^n \cap \rho_k B_2^n) \sim \sqrt{\log(en/k)}.
\end{equation}

Hence, combined with Lemma \ref{lemma:subgaussian-diameter} for $u$ and $\delta$ that depend only on $L$, one has that $Pr({\cal A}) \geq 1-2\exp(-c_4(L)k)$, and if $\sigma \in {\cal A}$ there is $I \subset \{1,...,N\}$, $|I| \geq c_5(L)k$, for which
$$
c_6(L) \sqrt{\frac{\log(en/k)}{k}} B_\infty^I \subset \Gamma_I (B_1^n \cap \rho_k B_2^n) \subset c_7(L) \sqrt{\frac{\log(en/k)}{k}} B_2^I.
$$

This observation should be compared with the following result from \cite{LPRT}:
\begin{Theorem} \label{thm:LPRT}
For every $L > 1$ there exist constants $c_1$ and $c_2$ that depend on $L$ and for which the following holds. Let $\xi$ be mean-zero, variance one, $L$-subgaussian random variable. Set $X=(\xi_i)_{i=1}^n$ to be a vector with independent coordinates, distributed according to $\xi$ and put $\Gamma=\sum_{i=1}^k \inr{X_i,\cdot}e_i$, where $(X_i)_{i=1}^k$ are independent copies of $X$. Then, for $0<\beta<1/2$, with probability at least $1-2\exp(-c_1k^\beta n^{1-\beta})$,
$$
c_2 \sqrt{\frac{\log(en/k)}{k}} B_\infty^k \subset \Gamma B_1^n .
$$
\end{Theorem}

\vskip0.5cm

\subsubsection*{Improving the lower estimate using a further projection}

Another outcome of Theorem \ref{thm:beyond-gaussian-main} is that
if $\sigma \in {\cal A}$, $Q_I(P_\sigma F)$ is only a proportional linear projection away from being equivalent to a Euclidean ball. This is well known for a typical orthogonal projection relative to the Haar measure, but we will show that the same is true for more general random ensembles.

Let $\sigma \in {\cal A}$, set $I$ to be as in Theorem \ref{thm:beyond-gaussian-main}, recall
that $|I|\sim k_F^*$ and put $W=Q_IV=Q_I(P_\sigma F)$.

For a mean-zero, variance 1, $L$-subgaussian random variable $\xi$, let $Y=(\xi_i)_{i=1}^{I}$. If $Y_1,...,Y_M$ are independent copies of $Y$,  consider what a further linear projection, $\sum_{i=1}^M \inr{Y_i,\cdot}e_i : \R^{I} \to \R^M$ does to both sides of the inclusion in \eqref{eq:beyond-gaussian-full}.

Since $Y$ is an isotropic, $L$-subgaussian random vector on $\R^{I}$, then by Lemma \ref{lemma:subgaussian-diameter}, with probability at least $1-2\exp(-c_0M)$,
$$
\sup_{v \in B_2^I} \left(\sum_{i=1}^M \inr{Y_i,v}^2 \right)^{1/2} \lesssim_{L} \sqrt{|I|} + \sqrt{M}.
$$
Hence, for $M \sim k_F^*$, and $\Gamma_\tau=\sum_{i=1}^M \inr{Y_i,\cdot}e_i$, one has $\Gamma_\tau W \subset c_1\sqrt{k_F^*} \E\|G\|_F B_2^M$.

For the reverse inclusion, it suffices to prove that $\Gamma_\tau (\E\|G\|_F/\sqrt{k_F^*}) B_\infty^I$ contains a large enough Euclidean ball.
\begin{Lemma} \label{lower-D-cube}
For every $L \geq 1$ there exist constants $c_1,c_2$ and $c_3$ that depend only on $L$ and for which the following holds.
Let $\xi$ be a mean-zero, variance $1$, $L$-subgaussian random variable, and let $Y=(\xi_i)_{i=1}^m$ be a vector of independent copies of $\xi$.
If $M=c_1 m$ and $\Gamma_\tau=\sum_{i=1}^M \inr{Y_i,\cdot}e_i$, then with probability at least $1-\exp(-c_2m)$,
$$
c_3 m B_2^M \subset \Gamma_\tau B_\infty^m.
$$
\end{Lemma}

\proof
Following the same path as in the proof of Theorem \ref{thm:lower-D-gaussian}, it suffices to show that $\inf_{s \in S^{M-1}} \|\sum_{i=1}^M s_i Y_i\|_{\ell_1^m} \geq \rho$, to ensure that $\rho B_2^M \subset \Gamma_\tau B_2^m$. Let $s \in S^{M-1}$, set $z=\sum_{i=1}^M s_i \xi_i$ and observe that if $z_1,...,z_N$ are independent copies of $z$ then $\sum_{i=1} s_i Y_i$ has the same distribution as $Z=(z_j)_{i=1}^m$.

It is standard to verify that since $z$ is mean-zero, variance one and $L$-subgaussian, $\E\|Z\|_{\ell_1^m} \sim m$. Applying Theorem \ref{thm:basic-probab-inequality},
$$
Pr\left( \frac{1}{m}\sum_{j=1}^m |z_j|  \leq \eps \right) \leq
Pr\left( \left|\frac{1}{m}\sum_{j=1}^m |z_j|-\E|z| \right|  \geq (\E|z|)/2 \right) \leq 2\exp(-c_1m),
$$
where $c_1$ depends only on $L$. The proof now follows from an $\eps$-net argument.

\begin{Corollary} \label{cor:full-D-general}
Using the notation above, if $N = c_1 k_F^*$, $\sigma \in A_{\delta,u,N}$, $I$ satisfies \eqref{eq:beyond-gaussian-full}, $W=Q_I(P_\sigma F)$ and $M=c_2k_F^*$, then with probability at least $1-2\exp(-c_3 M)$ with respect to $(Y_i)_{i=1}^{M}$,
$$
c_4 \sqrt{k_F^*} \E\|G\|_F B_2^{M} \subset \Gamma_\tau W \subset c_5 \sqrt{k_F^*} \E\|G\|_F B_2^{M}
$$
for constants that depend only on $L$.
\end{Corollary}

\vskip0.5cm

As an example, consider the case of a convex body $T \subset \R^n$ for which $\| \ \|_{T^\circ}$ has a nontrivial gaussian cotype 2 constant. If $|\tau| \sim k_T^*$ as in Corollary \ref{cor:full-D-general} and $W=\Gamma_I T$, it follows that
\begin{equation} \label{eq:dvo-cotype}
c_1\sqrt{k_T^*}\ell_*(T) B_2^{|\tau|} \subset \Gamma_\tau W \subset c_2\sqrt{k_T^*}\ell_*(T) B_2^{|\tau|},
\end{equation}
which is an isomorphic Dvoretzky-type theorem, obtained by first applying the random subgaussian operator $\Gamma_I=\sum_{i \in I} \inr{X_i,\cdot}e_i$ to $T$ for a specific choice of $I \subset \{1,...,N\}$ that depends on $(X_i)_{i=1}^N$, and is of cardinality proportional to $k_T^*$, and then a further linear projection given by $\sum_{i=1}^M \inr{Y_i,\cdot}e_i$, again, for $M$ that is proportional to $k_T^*$.

Equation \eqref{eq:dvo-cotype} should be compared with the Dvoretzky type theorem from  \cite{Men-TJ}. In \cite{Men-TJ} the random projection was given by a matrix with $N$ independent rows, generated by a random vector $X$ with iid mean-zero, variance $1$, $L$-subgaussian coordinates. The proof followed the path of a gaussian Dvoretzky Theorem -- with some modifications, and as such, it was based on a concentration argument. However, due to the lack of a strong enough concentration estimate, the dimension obtained in \cite{Men-TJ} was only $\sim_L k_T^*/\log^2 n$, and not $\sim_L k_T^*$ as we have here -- though the higher dimension comes at a price of an additional subgaussian projection.

To the best of our knowledge, it is still not known whether one may obtain an isomorphic Dvoretzky type theorem in such a case, using a single subgaussian projection generated by $X$ and of dimension $k_T^*$.

\footnotesize {
}


\begin{thebibliography}{10} \frenchspacing
%
\bibitem{Dud:book} R. M. Dudley, {\it Uniform Central Limit
Theorems}, Cambridge Studies in Advanced Mathematics 63, Cambridge
University Press, 1999.
%
\bibitem{F} X. Fernique, R\'{e}gularit\'{e} des trajectoires des fonctiones al\'{e}atoires
gaussiennes, Ecole d'Et\'{e} de Probabilit\'{e}s de St-Flour 1974,
Lecture Notes in Mathematics 480, 1-96, Springer-Verlag 1975.

\bibitem{Gor1} Y. Gordon, Some inequalities for Gaussian processes and applications, Israel Journl of Mathematics 50(4), 265-289, 1985.

\bibitem{Gor2} Y. Gordon, Gaussian processes and almost spherical sections of convex bodies, Annals of Probability 16(1), 180-188, 1988.


\bibitem{GLMP} Y. Gordon, A. Litvak, S. Mendelson, A. Pajor, Gaussian averages of interpolated bodies, Journal of Approximation Theory, 149, 59-73, 2008.

\bibitem{LT} M. Ledoux, M. Talagrand,  {\it Probability in
  Banach spaces. Isoperimetry and processes}, Ergebnisse der Mathematik
  und ihrer Grenzgebiete (3), vol. 23.  Springer-Verlag, Berlin, 1991.

\bibitem{LPRT} A. Litvak, A. Pajor, M. Rudelson, N. Tomczak-Jaegermann, Smallest singular value of random matrices and geometry of random polytopes, Advances in Mathematics, 195, 491-523, 2005.
%


\bibitem{MPR} S. Mendelson, A. Pajor, M. Rudelson, On the Geometry of random $\{-1,1\}$-polytopes, Discrete and Computational Geometry, 33(3) 365-379, 2005.

\bibitem{Men-TJ} S. Mendelson, N. Tomczak-Jaegermann, A subgaussian embedding theorem, Israel Journal of Mathematics, 164, 349-364, 2008.
%
\bibitem{Mil71} V.D. Milman, A new proof of A. Dvoretzky's theorem on cross-sections of convex bodies, Functional Analysis and its applications, 5(4), 28-37, 1971.
%
\bibitem{MilSch} V.D. Milman, G. Schechtman, { \it Asymptotic theory of
finite dimensional normed spaces}, Lecture Notes in Mathematics
1200, Springer, 1986.
%
%
\bibitem{Pis} G. Pisier {\it The Volume of Convex Bodies and Banach Space Geometry}, Cambridge Tracts in Mathematics vol 94, 1989.
%
%
\bibitem{RudVer} M. Rudelson, R. Vershynin, Combinatorics of random processes and sections of convex bodies, Annals of Mathematics 164, 603-648, 2006.
%
%
\bibitem{Tal87} M. Talagrand, Regularity of Gaussian processes,
Acta Math. 159, 99--149, 1987.
%
\bibitem{Tal:book} {M. Talagrand,} {\it upper and lower bounds for stochastic processes},
Springer, 2014.
\bibitem{VW}{A.W. Van der Vaart, J.A. Wellner, }
{\it Weak convergence and empirical processes}, Springer Verlag,
1996.
%

\end{thebibliography}
\end{document}